\documentclass[10pt]{article} 
\usepackage{amsfonts}
\usepackage{amssymb,amsfonts,amsmath,amsthm,cite,color}
\usepackage{enumerate}
\usepackage{array}

\PassOptionsToPackage{hyphens}{url}
\usepackage{hyperref}
\allowdisplaybreaks[4]

\newtheorem{thm}{Theorem}[section]
\newtheorem{cor}[thm]{Corollary}
\newtheorem{lem}[thm]{Lemma}
\newtheorem{pro}[thm]{Proposition}

\theoremstyle{definition}

\makeatletter \@addtoreset{equation}{section}

\def\pf{\noindent {\it Proof.\ }}

\title{\vspace{-1.6cm}
	On the Rank Functions of Powerful Sets}
\author{Benjamin R. Jones   \\
	Faculty of Information Technology   \\
	Monash University   \\
	Clayton, Victoria 3800   \\
	Australia   \\
	Email: \href{mailto:Benjamin.Jones@monash.edu}{\texttt{Benjamin.Jones@monash.edu}}}
\date{\today}

\begin{document}
	\maketitle
	
	\begin{abstract}

	A set $S\subseteq 2^E$ of subsets of a finite set $E$ is \emph{powerful} if, for all $X\subseteq E$, the number of subsets of $X$ in $S$ is a power of 2. 
	Each powerful set is associated with a non-negative integer valued function, 
	which we call the rank function. 
	Powerful sets were introduced by Farr and Wang as a generalisation of binary matroids, as the cocircuit space of a binary matroid gives a powerful set with the corresponding matroid rank function.
	
	In this paper we investigate how structural properties of a powerful set can be characterised in terms of
	its rank function. Powerful sets have four types of degenerate elements, including loops and coloops. We show that certain evaluations of the rank function of a powerful set determine the degenerate elements.
	We introduce powerful multisets and prove some fundamental results on them. 
	We show that a powerful set corresponds to a binary matroid if and only if its rank function is subcardinal. This paper answers the two conjectures made by Farr and Wang in the affirmative.

	\end{abstract}
	
	%-------------------------------------------------------------

	\section{Introduction}
	\label{sec:intro}
	Let $E$ be a finite set, and let $S\subseteq 2^E$. We say $S$ is \emph{powerful} if, for every $X\subseteq E$, the number of elements in $S$ that are subsets of $X$ is a power of 2.
	
	We call $E$ the \emph{ground set} of $S$, and define the \emph{rank function} of $S$ as the function $r_S:2^E\to\mathbb{R}$ defined by
	\[r_S(X)=\log_2\left(\frac{|S|}{|\left\{Y\in S~:~Y\subseteq E\backslash{X}\right\}|}\right).\]
	This rank function is integer valued if and only if $S$ is a powerful set, as the numerator and denominator are both powers of 2.
	
    Binary functions and their corresponding rank functions were introduced in \cite{farr93} as a generalisation of binary matroids. Powerful sets, introduced in \cite{farr2017}, are $\{0,1\}$-valued binary functions with integer valued rank functions.  We can think of a powerful set $S$ as a code by considering the characteristic vectors of sets in $S$. We will freely move between the viewpoints of a powerful set as a collection of subsets, and as a code.
		
  One particularly interesting class of powerful sets are the binary linear spaces, each of which is a powerful set. These \emph{linear} powerful sets can also be thought of as the cocircuit spaces of a binary matroid, with corresponding rank function. The \emph{nonlinear} powerful sets are those which do not correspond to binary matroids, and \cite{farr2017} showed that almost all powerful sets are nonlinear. 
  
    We extend the work started by \cite{farr2017}, focusing primarily on properties of the rank function of powerful sets. We express a number of structural properties of a powerful set in terms of the rank function.
    
    In Section \ref{sec:extensions}, we will consider a number of the special elements introduced by \cite{farr2017} that extend powerful sets, and determine how these operations are expressed in terms of the rank function. One particular type of element, the star, is the focus of Section \ref{sec:star}.

    Contraction and deletion operations on powerful sets were also defined in \cite{farr2017}.
    One issue with the deletion operation is that deletion from a powerful set does not necessarily return a powerful set. Section \ref{sec:multisets} introduces powerful multisets as a generalisation of powerful sets where the deletion operation is well behaved. We prove some elementary results on powerful multisets, and use these objects to determine properties of powerful sets in the later sections. Section \ref{sec:nonlinear} proves that a powerful set is linear if and only if its rank function is subcardinal.
    
    We refer the reader to \cite{oxley92,welsh76} for an overview of matroid theory.
    
	\section{Degenerate elements}\label{sec:extensions}
	
  Loops and coloops are the ``degenerate'' elements in matroid theory. Farr and Wang \cite{farr2017} identified four types of degenerate element in powerful sets, two of which are analogous to loops and coloops. Here, we will define these four types and state some elementary results.

Let $S\subseteq 2^E$, and $e\in E$. We say $e$ is a \emph{loop} if $e$ is in no member of $S$. There will then be a set $T\subseteq 2^{E\backslash \{e\}}$ such that $S=\{X:X\in T\}$, and we write $S=T+\circ_e$.
	
We say $e$ is a \emph{coloop} if there exists $T\subseteq 2^{E\backslash\{e\}}$ such that $S=\{X,X\cup\{e\}:X\in T\}$, and we write $S=T+\circ^*_e$.
	
	We say $e$ is a \emph{frame} if there exists $T\subseteq 2^{E\backslash\{e\}}$ such that $S=\{X\cup\{e\}:X\in T\backslash\{\emptyset\}\}\cup\{\emptyset\}$, and we write $S=T+\Box_e$.
	
	We say $e$ is a \emph{star} if there exists $T\subseteq 2^{E\backslash\{e\}}$ such that $S=\{X:X\in T\}\cup\{X\cup\{e\}:X\notin T\}$, and we write $S=T+\star_e$.
		
	\begin{pro}[Theorems 3.1, 3.2, 3.5, 3.7, \cite{farr2017}]
	  Let $T\subseteq 2^{E\backslash\{e\}}$. Then
	  
	\begin{enumerate}[a)] 
	    \item $T+\circ_e$\, is powerful if and only if $T$ is powerful,
	    
	    \item $T+\circ^*_e$\, is powerful if and only if $T$ is powerful,

	    \item $T+\Box_e$ is powerful if and only if $T$ is powerful,

	    \item $T+\star_e$\; is powerful if and only if $T$ is powerful.
	\end{enumerate}
	\end{pro}
	
\subsection{Results on direct sum}

	We have the following operation on pairs of powerful sets from \cite{farr2017}.  Given powerful sets $S$ and $T$, with disjoint ground sets $E$ and $F$ respectively, the \emph{direct sum}  $S\oplus T$, on ground set $E\cup F$ is given by
\[		\{X\cup Y\subseteq E\cup F~:~X\in S,~ Y\in T	\}.	\]
Loops and coloops are special cases of the direct sum of powerful sets. A powerful set $T+\circ_e$ on ground set $E$ with loop $e$ is the direct sum $T\oplus \{\emptyset\}$ of powerful sets $T$ and $\{\emptyset\}$ with ground sets $E\backslash \{e\}$ and $\{e\}$ respectively. Likewise for a coloop, $T+\circ^*_e= T \oplus \{\emptyset,\{e\}\}$.

	The minimal nonempty members of a powerful set $S$ are called the \emph{cocircuits} of $S$, and we denote the set of cocircuits of $S$ by $\mathcal{C}(S)$.
	In \cite{farr2017}, the following theorem was proved.
	\begin{thm}\label{thm:cocircuits}
		Every powerful set is determined by its set of cocircuits.
	\end{thm}
	We use this in the proof of Theorem \ref{thm:coloops}, and later in \S \ref{sec:nonlinear}.

	The following theorem relates to direct sums of powerful sets, and suggests that powerful sets have structures analogous to components in binary matroids. Loops and coloops are specialisations of the direct sum, and so observations on the direct sum will serve us well. %hmm
	There are some immediate consequences of the following theorem, including the proof of Conjecture 5 in \cite{farr2017}.
	\begin{thm} \label{thm:direct_sum}Let $S$ and $T$ be powerful sets with non-intersecting ground sets. Then,
		$\mathcal{C}(S\oplus T)=\mathcal{C}(S)\cup\mathcal{C}(T)$.\end{thm}
	
	\pf
	$(\subseteq)$ Take some nonempty $X\cup Y\in S\oplus T$, with $X\in S,Y\in T$.

	Suppose $X$ is nonempty and not a cocircuit of $S$. Then there is some $C\in\mathcal{C}(S)\subseteq S\oplus T$ such that $C\subset X\subseteq X\cup Y$, and so $X\cup Y$ is not a cocircuit. Otherwise, $Y$ is nonempty, and repeating this argument gives the cocircuits of $S$ and $T$ as the only potential cocircuits of $S\oplus T$.
	
	Hence $\mathcal{C}(S\oplus T)\subseteq \mathcal{C}(S)\cup\mathcal{C}(T)$.
	
	($\supseteq$) Any member of $\mathcal{C}(S)\cup\mathcal{C}(T)$ is a minimal nonempty member of $S\oplus T$, and so $\mathcal{C}(S\oplus T)\supseteq \mathcal{C}(S)\cup\mathcal{C}(T) $. $\hfill\Box$
	
\subsection{Characterisation of coloops, loops and frames}
	For binary matroids, loops and coloops are the elements $e$ where $r(\{e\})=0$ and $r(E)-r(E\backslash \{e\})=1$ respectively. These degenerate elements are thus identifiable by their rank. This is also the case for loops, coloops and frames in the more general setting of powerful sets. This section will prove this, resolving a conjecture from \cite{farr2017}.

	\begin{thm}\label{thm:coloops}
		An element $e\in E$ is a coloop of a powerful set $S\subseteq 2^E$ if and only if $\{e\}\in S$.	
	\end{thm}
\pf $(\Rightarrow)$ Suppose $e$ is a coloop of $S$, then we write $S=T+\circ_e^*$. As $T$ is a powerful set, 
$\emptyset\in T$ and so $\{e\}\in S$.

$(\Leftarrow$) Suppose $\{e\}\in S$. Let $T\subseteq 2^{E\backslash\{e\}}$ be the set of elements in $S$ that do not contain $e$. We have that $T$ is a powerful set, as for any $X\subseteq E\backslash\{e\}$ the number of elements in $T$ that have trivial intersection with $X$ is equal to the number of elements in $S$ that have trivial intersection with $X\cup\{e\}$, which is a power of two since $S$ is powerful.

We have $\mathcal{C}(S)=\mathcal{C}(T)\cup\{\{e\}\}=\mathcal{C}(T)\cup\mathcal{C}(\{\emptyset,\{e\}\})$ and so by Theorems \ref{thm:cocircuits} and \ref{thm:direct_sum} we have $S=T+\circ^*_e$. Hence $e$ is a coloop of $S$.    $\hfill\Box$\\

This theorem answers Conjecture 5 from \cite{farr2017} in the affirmative. 
	
	We can characterise these degenerate elements in terms of evaluations of the rank function of a powerful set.
	\begin{thm} 
		\label{thm:extension_rank}
		Let $S\subseteq 2^E$ be a powerful set, and $e\in E$. Then,
		\begin{enumerate}
			\item The element $e$ is a loop if and only if $r(\{e\})=0$,
			\item The element $e$ is a frame if and only if $r(\{e\})=r(E)$,
			\item The element $e$ is a coloop if and only if $r(E)-r(E\backslash \{e\})=1$.
		\end{enumerate}
	\end{thm}
	\pf {
	Points \emph{1.} and \emph{2.} follow immediately from the definition of loops and frames. Point \emph{3.} follows from Theorem \ref{thm:coloops}}. $\hfill\Box$ 
	
	Stars are the only type of degenerate element mentioned that have not been covered in this manner, which will be addressed in Section \ref{sec:star}. But first we will introduce powerful multisets.
	\section{Powerful multisets}\label{sec:multisets}
	In this section, we extend the definition of powerful sets to powerful multisets and establish some properties that will be used in Sections \ref{sec:star} and \ref{sec:nonlinear}. Powerful multisets are worthy of investigation in their own right, but also serve as a tool to study properties of powerful sets.

We will consider a multiset, over finite ground set $E$, to be a collection of subsets of $E$ with repetition. For such a multiset $S$, we can define an \emph{indicator function} $f:2^E\to\mathbb{N}\cup\{0\}$, where $f(X)$ is the number of times $X$ is repeated in $S$, for each $X\subseteq E$.

A multiset $S$, with ground set $E$ and indicator function $f$, is \emph{powerful} if, for all $X\subseteq E$, the rank of $X$, given by
	\[  	Qf(X)=\log_2\left(\frac{\sum\limits_{Y\subseteq E}f(Y)}{\sum\limits_{Y\subseteq E\setminus X}f(Y)}\right),\]
	is an integer. The transform $Q$ comes from \cite{farr93}.
	
	Note for the indicator function $f$ that $f(\emptyset)\neq 0$ if and only if $Qf$ is a well defined function. From now on we will assume that a given multiset contains at least one copy of the empty set when we calculate its rank function. 

Two powerful multisets differing by a nonzero scalar will have the same rank function, as $Q(\alpha f)=Qf$ for $\alpha\in \mathbb{R}\backslash\{0\}$.	
 We call two powerful multisets $S_1$ and $S_2$, with ground sets $E_1$ and $E_2$ and indicator functions $f_1$ and $f_2$ respectively, \emph{isomorphic} if there is a bijective function $\phi:E_1\to E_2$ and constant $\alpha\in \mathbb{R}\backslash\{0\}$ such that, for all $X\subseteq E_1$,  $f_1(X)=\alpha f_2(\phi(X))$.

  We use the same notation for powerful multisets as for powerful sets. This includes writing the rank function $Qf$ as $r_S$. We will make clear whether $S$ is a set or multiset.

	The empty set is always an element of a powerful set, and this also holds for powerful multisets. However, the multiplicity of the empty set may be greater than one for a powerful multiset. We will now show that each powerful multiset is isomorphic to one where the empty set has multiplicity one.	
	
	\begin{thm}
			\label{thm:multiset_scale}
		If $S$ is a powerful multiset with indicator function $f$, then, for all $X\subseteq E$, $f(\emptyset)$ divides $f(X)$.
	\end{thm}
\pf{
First, we note that, for any $X\subseteq E$,
\[		r_S(E\backslash X)=\log_2\left(\frac{\sum\limits_{Y\subseteq E}f(Y)}{\sum\limits_{Y\subseteq X}f(Y)}\right)=r_S(E)-\log_2\left(\frac{\sum\limits_{Y\subseteq X}f(Y)}{f(\emptyset)}\right)\in \mathbb{N}\cup\{0\}.\]
Since $r_S(E),r_S(E\backslash X) \in \mathbb{N}\cup\{0\}$ and $r_S(E\backslash X)\leq r_S(E)$, we have $f(\emptyset) \mid \sum\limits_{Y\subseteq X}f(Y)$.

Suppose for contradiction that the theorem is not true, and there is some minimal $X\subseteq E$ where $f(\emptyset) \nmid f(X)$. We have $ \sum\limits_{Y\subseteq X}f(Y) =f(X)+\sum\limits_{Y\subset X}f(Y)$, and so
\[f(\emptyset) \mid \left(f(X)+\sum\limits_{Y\subset X}f(Y)\right).\]

Since $X$ is minimal, $f(\emptyset) \mid f(Y)$ for all $Y\subset X$, and so $f(\emptyset) \mid \sum\limits_{Y\subset X}f(Y)$. Hence $f(\emptyset) \mid f(X)$, which contradicts our assumption.$\hfill\Box$}
	
Given a powerful multiset $S$ with indicator function $f$, this theorem guarantees that $\frac{1}{f(\emptyset)}f$ is the indicator function of an isomorphic powerful multiset (noting that $f(\emptyset)$ is nonzero, as $Qf(E)$ is well defined).

\begin{pro}
    A powerful multiset $S$, with indicator function $f$, is isomorphic to a powerful set if and only if $\frac{1}{f(\emptyset)}f$ is $\{0,1\}$-valued. 
\end{pro}

The deletion and contraction operations of matroids were extended to multisets in \cite{farr93}. Let $S$ be a powerful multiset with ground set $E$ and indicator function $f$, and take $e\in E$. We define the function $f/ e$ where, for each $X\subseteq E\backslash\{e\}$,
\[  f/e(X)=f(X).     \]
The multiset $S/e$, with indicator function $f/e$, is formed by the \emph{contraction} of $e$ from $S$. Similarly, we define $f\backslash e$ where, for each $X\subseteq E\backslash\{e\}$,
\[  f\backslash e(X)={f(X)+f(X\cup \{e\})},\]
and the corresponding multiset $S\backslash e$ as being formed by the \emph{deletion} of $e$ from $S$.

\begin{pro}[Theorem 4.5 \cite{farr93}]\label{pro:minor_rank} For a multiset $S$ with ground set $E$, for all $e\in E$ and $X\subseteq E\backslash\{e\}$
\[\begin{array}{ll}
     r_{S/e}(X)=&r_S(X\cup\{e\})-r_S(\{e\}),  \\
     r_{S\backslash{e }} (X)=&r_S(X).    \\
\end{array}\]
\end{pro}

These deletion and contraction operations for powerful multiset rank functions are the same as those for matroid rank functions. Extensions of these operations to ``arbitrary'' rank functions have also been considered in \cite{gordon2012}. 

\begin{cor}
    Let $S$ be a powerful set with ground set $E$. For every $e\in E$, both $S/e$ and $S\backslash e$ are powerful multisets.
\end{cor}

It is important to look at how the contraction and deletion operations work on powerful sets.

For a powerful set $S\subseteq 2^E$, contraction \cite{farr2017} of $e\in E$ from $S$ is the powerful set given by 
\[      S/e=\{X\subseteq E\backslash\{e\}:X\in S\text\}.  \]
If deletion of $e$ from $S$ gives a powerful set, then $S\backslash e$ is of the form
\[  S\backslash e =\{ X\subseteq E\backslash\{e\}:X\in S \text{ or } X\cup \{e\}\in S    \}.      \]  
However, deletion from a powerful set does not always give a powerful set. For example, for the powerful set $S=\{\emptyset,\{1,3\},\{2,3\},\{1,2,3\}\}$, we have that $S\backslash 1$ is the powerful multiset $S\backslash 1 =\{\emptyset,\{3\},\{2,3\},\{2,3\}\}$ which is not a set.  

For a powerful set $S\subseteq 2^E$, we call an element $e\in E$  \emph{deletable} if $S\backslash{e}$ is a powerful set. Equivalently, where $f$ is the $\{0,1\}$-valued indicator function of $S$, an element $e$ is deletable if  $1+f(\{e\})$  divides $f(X)+f(X\cup\{e\})$, for every $X\subseteq E\backslash\{e\}$.

If $S$ is linear, then every element is deletable. The converse, however, does not hold, as every element of the nonlinear powerful set $S=\left\{\emptyset,\{1,2,4\},\{1,3,4\},\{2,3,4\}\right\}$ is deletable.

Deletion and contraction in powerful sets can affect the size of the powerful set as follows.

	\begin{pro}\label{pro:minor_size} For a powerful set $S\subseteq 2^E$ and $e\in E$, we have
\[  |S/e|=|S|/2^{r_S(\{e\})}.\]
If $e$ is deletable, then
\[ |S\backslash e|=
\left\{\begin{array}{rl}
   \frac{1}{2}|S|  &\text{if } e\text{ is a coloop}, \vspace{1.5mm}\\
   |S| &    otherwise.
\end{array}         \right.\]
    
\end{pro}
\pf{ Let $f$ be the $\{0,1\}$-valued indicator function of the powerful set $S$. The indicator function of $f/e$, given by $f/e(X)=f(X)$, is $\{0,1\}$-valued for each $X\subseteq E\backslash e$. Hence 
\[  |S/e|=\sum\limits_{X\subseteq E\backslash e} f/e(X) = \left(\frac{|S|}{\sum\limits_{X\subseteq E} f(X)} \right) \sum\limits_{X\subseteq E\backslash e} f(X)  =|S|/ 2^{r_S(\{e\})}. \]
If $e$ is deletable, then the indicator function $f\backslash e$ of the powerful set $S\backslash e$ is given by $f\backslash e(X)=f(X)+f(X\cup \{e\})$. Note that $f\backslash e$ is not necessarily $\{0,1\}$-valued, but $\frac{1}{f\backslash e(\emptyset)}f$ is. 
Hence
\[ |S\backslash e|= \frac{1}{f\backslash e(\emptyset)} \sum\limits_{X\subseteq E\backslash e} f\backslash e(X)= \frac{1}{1+f(\{e\})} \sum\limits_{X\subseteq E\backslash e}\left(f(X)+f(X\cup \{e\})\right)=\frac{|S|}{1+f(\{e\})}.    \]
 By Theorem \ref{thm:coloops}, $f(\{e\})$ is 1 if $e$ is a coloop, and zero otherwise.     $\hfill\Box$}
	
	\section{Stars}\label{sec:star}
In this section we give a characterisation for the star elements of powerful sets.   Lemma  \ref{starcols}  is used in Theorem \ref{thm:deletable_rank} to prove the existence of deletable elements of some powerful sets.
	\begin{lem}
		\label{starcols}
		Let $S$ be a powerful set of order $n$ and rank $n-1$. There is at most one element $e\in E$ which has rank not equal to 1.
	\end{lem}
\pf

 Suppose that for some $e\in E$ we have $r(\{e\})=0$. By Theorem \ref{thm:extension_rank}, $e$ is a loop, and so $S$ must be isomorphic to $T+\circ_e$ for some powerful set $T\subseteq 2^{E\backslash\{e\}}$.  By Proposition \ref{pro:minor_rank}, $r_T(E\backslash\{e\})=r_{{T+\circ_e}}(E)-r_{{T+\circ_e}}(\{e\})=n-1-0$, and so $T$ must be equal to $2^{E\backslash\{e\}}$. Note that for each $X\subseteq E\backslash \{e\}$, $r_{2^{E\backslash\{e\}}}(X)=|X|$.
 
For each $e'\in E\backslash\{e\}$,   $r_S(\{e'\})=r_{2^{E\backslash\{ e\}}+\circ_e}(\{e'\})=r_{2^{E\backslash\{ e\}}}(\{e'\})=1$. Hence, all elements of $E$ except $e$ have rank $1$. 

Suppose then that $S$ has no element of rank 0. By Theorem \ref{thm:extension_rank}, the rank of each element of $E$ is at least 1.  Suppose now for contradiction that there is a powerful set $S$ of order $n\geq 2$ and rank $n-1$ with distinct elements $e,f \in E$ such that $r(\{e\}),r(\{f\})\geq 2$. 
 
 We have
 
 \begin{center}

 $ \begin{array}{llll}
|S| &   =& |\{X\subseteq E\backslash\{e,f\}:X\cup\{e,f\}\in S\}|+|\{X\subseteq E\backslash\{e\}:X\in S\}|\\  &&~~+|\{X\subseteq E\backslash\{f\}:X\in S\}|-|\{X\subseteq E\backslash\{e,f\}:X\in S\}|              \\
    &   =& |\{X\subseteq E\backslash\{e,f\}:X\cup\{e,f\}\in S\}| +|S/e|+|S/f|-|S/e/f|    \\
    &   \leq& |2^{E\backslash\{e,f\}}| +|S/e|+|S/f|-|S/e/f|    \\
    &   = &|2^{E\backslash\{e,f\}}| +2^{r(E)-r(\{e\})}+2^{r(E)-r(\{f\})}-|S/e/f|    \\
    & &~\hbox{(by Proposition \ref{pro:minor_size})}\\
    &   \leq&  2^{n-2}+2^{n-3}+2^{n-3}-1\\
    &  &~\hbox{(since $r(E)=n-1$, $r(\{e\}),r(\{f\})\geq 2$, and $\emptyset\in S/e/f$)}\\
 & =& 2^{n-1}-1.
\end{array}$
\end{center}

But $|S|=2^{n-1}$, which is a contradiction. $\hfill\Box$
 
\vspace{5mm}

We introduce some notation which we will use in the proof of Theorem \ref{thm:deletable_rank}. Let $T$ be a multiset on ground set $E$ with indicator function $f_T$. For $X\subseteq E$, let 
\begin{equation} 
    z_X(T)= \sum\limits_{Y\subseteq {E\backslash{ X}}}f_T(Y).
\end{equation}
For a set $T$, $z_X(T)$ counts the number of sets $Y\in T$ where $|X \cap Y|=0$. The following proposition presents some basic properties of set systems using this notation.  
\begin{pro}~\label{pro:z_stuff}
    Let $S,T\subseteq 2^E$.
\begin{enumerate}[ a)] 
    \item  For any $X\subseteq E$, $r_S(X)=\log_2\left(  z_\emptyset(S)/z_{X}(S)\right)$.
    
    \item   $S$ is powerful if and only if $z_X(S)$ is a power for two for every $X\subseteq E$.
    \item If, for every $X\subseteq E$, $z_X(S)=z_X(T)$, then $S=T$.

    \item  If, for every $X\subseteq E$, $z_X(S)+z_X(T)= 2^{|E|-|X|}$, then $S=2^E\backslash T$.
    
\end{enumerate}
\end{pro}

Note that if $e$ is a deletable element of a powerful set $S$ of order $n$ and rank $n-1$, and $e$ is not a coloop, then $S\backslash e$ is the powerful set $2^{E\backslash\{e\}}$. 

 	\begin{thm} \label{thm:deletable_rank}
 		Let $S\subseteq 2^E$ be a powerful set of order $n$ and rank $n-1$. Then there exists a deletable element $e\in E$ of $S$.
 	\end{thm}

 \pf

For contradiction, assume there exists a powerful set $S$ of order $n$ and rank $n-1$ that has no deletable elements. Take $S$, with ground set $E$, to be of minimal order.

Any powerful set of order $n\leq 2$ and rank $n-1$ is linear,  so every element is deletable. Therefore $S$ has order at least 3.

We first remark that as $S$ has no deletable elements,  $S$ must not contain any coloops, loops or frames. By Theorem \ref{thm:extension_rank} each element in $E$ has rank in $S$ not equal to $0$ or $n-1$. We now consider two cases based on the rank of single elements of $E$ in $S$.

~\\\noindent\textbf{Case 1:} There is an element $e\in E$ such that $r_S(\{e\})\geq 2$.

By Lemma \ref{starcols}, every other element of $E$ has rank 1.

\textbf{Case 1a:} There exist distinct $f_1,f_2\in E\backslash \{e\}$ and sets $Y_1,Y_2 \in S/e$ such that $f_1\in Y_1$ and  $f_2\in Y_2$.

We have $r_S(\{e\})=\log_2(|S|/z_{\{e\}}(S))$, so $z_{\{e\}}(S)=2^{n-1-r_S(\{e\})}$.

Consider now, $r_{S/f_1}(\{e\})= \log_2(z_\emptyset(S/f_1)/z_{\{e\}}(S/f_1))$. We have $z_\emptyset (S/f_1)=|S/f_1|=2^{n-2}$, as $f_1$ has rank 1 in $S$ (by Proposition \ref{pro:minor_size}). We also have $z_{\{e\}}(S/f_1)<z_{\{e\}}(S)=2^{n-1-r_S(\{e\})}$, as $Y_1\notin S/f_1$ and $Y_1\in S$. Therefore, $r_{S/f_1}(\{e\})>r_S(\{e\})-1\geq 1$.

We have $S/f_1$ is a powerful set of order $n-1$ and rank $r_{S/f_1}(E\backslash\{e\})=r_S(E)-r_S(\{f_1\})=(n-1)-1=n-2$ by Proposition \ref{pro:minor_rank}. By the minimality of the order of $S$, there exists a deletable element $x$ of $S/f_1$.

As $x$ is not a coloop of $S$, $x$ is not a coloop of $S/f_1$  (by Theorem \ref{thm:coloops}). So by Proposition \ref{pro:minor_size}, $|S/f_1\backslash x|=|S/f_1|=|S|/ 2^{r_S(\{f_1\})}=2^{n-2}$, and so  $(S/f_1)\backslash x =2^{E\backslash\{f_1,x\}}$.

We must have $x=e$, as otherwise $r_{(S/f_1)\backslash x}(\{e\})=r_{(S/f_1)}(\{e\})\neq 1$.
Therefore, $e$ is a deletable element of $S/f_1$, and $S/f_1\backslash{e} = 2^{E\backslash{\{e,f_1\}}}$. 

Similarly, $e$ is a deletable element of $S/f_2$, and  $S/f_2\backslash{e} = 2^{E\backslash{\{e,f_2\}}}$.

There are eight possible intersections with $\{e,f_1,f_2\}$ that a member of $S$ can have. We partition $S$ according to these intersections, giving the following:
 \begin{center}
 $
 	\begin{array} {c|c}
 	    Y       &     \{X\subseteq E\backslash\{e,f_1,f_2\}\ |\ Y\cup X \in S\}    \\\hline
 	\{\phantom{e}\phantom{,}\phantom{f_1}\phantom{,}\phantom{f_2}\} & T_1	\\
 	\{\phantom{e,}\phantom{f_1,}{f_2}\} & T_2   \\
 	\{\phantom{e,}{f_1}\phantom{,}\phantom{f_2}\} & T_3	\\	
 	\{\phantom{e,}{f_1,}{f_2}\}         & V\\	
 	\{{e}\phantom{,}\phantom{f_1}\phantom{,}\phantom{f_2}\}   & 2^{E\backslash\{e,f_1,f_2\}}\backslash T_1	\\
 	\{{e}{,}\phantom{f_1}\phantom{,}{f_2}\}         & 2^{E\backslash\{e,f_1,f_2\}}\backslash T_2	\\
 	\{{e}{,}{f_1}\phantom{,}\phantom{f_2}\}           & 2^{E\backslash\{e,f_1,f_2\}}\backslash T_3	\\
 	\{{e}{,}{f_1}{,}{f_2}\}                    & W	\\
 	\end{array}
 	$
\end{center}

Each row corresponds to a subset $Y\subseteq \{e,f_1,f_2\}$ on the left hand side, with the right hand side being the set of elements $X\subseteq E\backslash\{e,f_1,f_2\}$ such that $Y\cup X$ is in $S$.

Here, $T_1,T_2,T_3,V,W$ are subsets of $2^{E\backslash\{e,f_1,f_2\}}$. Note that the expressions in the fifth, sixth and seventh rows are derived from the facts that $S/f_1\backslash{e} = 2^{E\backslash{\{e,f_1\}}}$ and  $S/f_2\backslash{e} = 2^{E\backslash{\{e,f_2\}}}$.

 For any $X\subseteq  E\backslash \{e,f_1,f_2\}$,
 \[z_X(S)=3z_X(2^{E\backslash\{e,f_1,f_2\}})+z_X(V)+z_X(W)=3(2^{n-3-|X|})+z_X(V)+z_X(W).\]
 
  As $V,W\subseteq 2^{E\backslash\{e,f_1,f_2\}}$, $0\leq z_X(V)+z_X(W)\leq  2(2^{n-3-|X|})$, and as $z_X(S)$  is a power of 2 we must have $z_X(V)+z_X(W)=2^{n-3-|X|}$ for each $X\subseteq\{e,f_1,f_2\}$. By Proposition \ref{pro:z_stuff}{d}, we must have $W=2^{E\backslash\{e,f_1,f_2\}}\backslash V$, and so  $S\backslash e=2^{E\backslash\{e\}}$. Hence $e$ is a deletable element of $S$. 
  
  \textbf{Case 1b:} For all distinct $f_1,f_2\in E\backslash\{e\}$, there does not exist $Y_1,Y_2\in S/e$ such that $f_1\in Y_1$, and $f_2\in Y_2$.
  
  It follows that every member of $S/{e}$ has size at most 1. 
  
 As $e$ is not a deletable element of $S$, $e$ is not a frame, and so $S/e$ contains at least two elements, one of which must be a singleton set, say $\{g\}\in S/{e}$. Hence we have $\{g\}\in S$. By Theorem \ref{thm:coloops}, $g$ is a coloop of $S$, and so $S= S/g +\circ_g^*$. The powerful set $S/g$ has order $n-1$ and rank $n-2$, and so contains a deletable element $h\in E\backslash\{g\}$ by the minimality of $S$. So $S/g\backslash{h}$ has order and rank $n-2$, so is equal to $2^{E\backslash\{g,h\}}$. We have
 
\[
 \begin{array}{lllr}
  2^{E\backslash\{h\}}&=&2^{E\backslash\{g,h\}}+\circ_g^*&\text{as every element of the powerful set $2^{E\backslash\{h\}}$ is a coloop,}\\
  &=& S/g\backslash{h}+\circ_g^* & \text{as $S/g\backslash{h}=2^{E\backslash\{g,h\}}$,}\\
  &=& (S/g+\circ_g^*)\backslash{h} &\text{as deletion and coloops commute,}\\
  &=& S\backslash{h} &\text{ as $g$ is a coloop of $S$.}
 \end{array}
 \]
Therefore $h$ is a deletable element of $S$.\\~\\ 
 \noindent\textbf{Case 2:} For all $e\in E$, $r_S(\{e\})=1$.
 
We must have that $S$ has no coloops, as any coloop is a deletable element.
 
For each $e\in E$, $S/e$ is a powerful set of order $n-1$ and rank ${n-2}$ (by Proposition \ref{pro:minor_rank}), so by the minimality of $S$ there is some deletable element of $S/e$. Choose an $e\in E$ and let $f_1\in E\backslash \{e\}$ be a deletable element of $S/e$. Let $f_2\in E\backslash \{f_1\}$ be a deletable element of $S/f_1$. As $f_1$ and $f_2$ are not coloops of $S$, deletion of coloops does not change the rank (Proposition \ref{pro:minor_rank} and deletion and contraction commute,  both $S/e\backslash{f_1}$ and $S/f_1\backslash{f_2}$ have rank $n-2$. 
 
 \textbf{Case 2a:} $f_1$ and $f_2$ are distinct.
 
There are eight possible intersections with $\{e,f_1,f_2\}$ that a member of $S$ can have. We can then partition $S$ according to these intersections, giving the following:

 \begin{center}	
  	$
 	\begin{array} {c|c}
 	    Y       &     \{X\subseteq E\backslash\{e,f_1,f_2\}\ |\ Y\cup X \in S\}    \\\hline
 	\{\phantom{e}\phantom{,}\phantom{f_1}\phantom{,}\phantom{f_2}\} & T_1	\\
 	\{\phantom{e,}\phantom{f_1,}{f_2}\} & T_2   \\
 	\{\phantom{e,}{f_1}\phantom{,}\phantom{f_2}\} & 2^{E\backslash\{e,f_1,f_2\}}\backslash T_1	\\	
 	\{\phantom{e,}{f_1,}{f_2}\}         & 2^{E\backslash\{e,f_1,f_2\}}\backslash T_2	\\	
 	\{{e}\phantom{,}\phantom{f_1}\phantom{,}\phantom{f_2}\}   & T_3	\\
 	\{{e}{,}\phantom{f_1}\phantom{,}{f_2}\}         & 2^{E\backslash\{e,f_1,f_2\}}\backslash T_3	\\
 	\{{e}{,}{f_1}\phantom{,}\phantom{f_2}\}           & V	\\
 	\{{e}{,}{f_1}{,}{f_2}\}                    & W	\\
 	\end{array}
 	$
 \end{center}

Here, $T_1,T_2,T_3,V,W$ are subsets of $2^{E\backslash\{e,f_1,f_2\}}$. Note that the expressions in the third, fourth and sixth rows are derived by the facts that $S/e\backslash f_1 =2^{E\backslash\{e,f_1\}}$ and $S/f_1\backslash f_2 =2^{E\backslash\{f_1,f_2\}}$.

 For all $X\subseteq E\backslash\{e,f_1,f_2\}$, we have
 \begin{center}
 	$
 	\begin{array} {rl}
 	z_X(S)	=	&	\sum_{i=1}^3 \left(z_X(T_i)+z_X(2^{E\backslash\{e,f_1,f_2\}}\backslash T_i) \right)+z_X(V)+z_X(W)\\
 	        =   &   3z_X(2^{E\backslash\{e,f_1,f_2\}})+z_X(V)+z_X(W)\\
        	=   &  3(2^{n-3-|X|})+z_X(V)+z_X(W)			\\
 	\end{array}
 	$
 \end{center}
 As $V$ and $W$ are subsets of $2^{E\backslash\{e,f_1,f_2\}}$ we have $0 \leq z_X(V),z_{X}(W) \leq 2^{n-3-|X|}$, and so
 \[ 3(2^{n-3-|X|}) \leq z_X(S)\leq  3(2^{n-3-|X|})+2(2^{n-3-|X|} ).\]

 As $S$ is a powerful set, $z_X(S)$ is a power of two, so $z_X(S)=2^{n-1-|X|}$. Hence $z_X(V)+z_{X}(W)=2^{n-3-|X|}=z_X(2^{E\backslash\{e,f_1,f_2\}})$ for all $X\subseteq E\backslash\{e,f_1,f_2\}$, and so $W=2^{E\backslash\{e,f_1,f_2\}}\backslash V$, by Proposition \ref{pro:z_stuff}.
 
 Therefore $f_1$ is a deletable element of $S$.

 \textbf{Case 2b:} $f_2=e$.
 
We have that $S/f_1\backslash e = S/e\backslash f_1 = 2^{E\backslash\{e,f_1\}}$. Take some arbitrary $f_3\in E\backslash\{e,f_1\}$. We partition $S$ according to the intersections of members of $S$ with $\{e,f_1,f_3\}$, giving the following:	
 \begin{center}
    $
 	\begin{array} {c|c}
 	Y & \{X\subseteq E\backslash\{e,f_1,f_3\}\ |\ Y\cup X \in S\}   \\
 	\hline
 	\{\phantom{e}\phantom{,}\phantom{f_1}\phantom{,}\phantom{f_3}\} & T_1	\\
 	\{\phantom{e,}\phantom{f_1,}{f_3}\} & T_2   \\
 	\{\phantom{e,}{f_1}\phantom{,}\phantom{f_3}\} & 2^{E\backslash\{e,f_1,f_3\}}\backslash T_1	\\	
 	\{\phantom{e,}{f_1,}{f_3}\}         & 2^{E\backslash\{e,f_1,f_3\}}\backslash T_2	\\	
 	\{{e}\phantom{,}\phantom{f_1}\phantom{,}\phantom{f_3}\}   & 2^{E\backslash\{e,f_1,f_3\}}\backslash T_1	\\	
 	\{{e}{,}\phantom{f_1}\phantom{,}{f_3}\}         & 2^{E\backslash\{e,f_1,f_3\}}\backslash T_2	\\
 	\{{e}{,}{f_1}\phantom{,}\phantom{f_3}\}           & V	\\
 	\{{e}{,}{f_1}{,}{f_3}\}                    & W	\\
 	\end{array}
 	$
\end{center}
where $T_1, T_2,V,W$ are subsets of $2^{E\backslash\{e,f_1,f_3\}}$.
For $X\subseteq E\backslash\{e,f_1,f_3\}$, \begin{center}
 	$
 	\begin{array} {rl}
 	z_{X\cup\{f_3\}}(S)	=	&	z_X({T_1})+2z_X(2^{E\backslash\{e,f_1,f_3\}}\backslash T_1)+z_X(V)			\\
 	=	&	2^{n-2-|X|}-z_X({T_1})+z_X(V).\\
 	\end{array}
 	$
 \end{center}
 Since $0\leq z_X(V),z_X({T_1})\leq 2^{n-3-|X|}$ and $z_{X\cup\{f_3\}}(S)$ is a power of 2, we have $z_X(V)=z_{X}({T_1})$ for all $X$, and so $V=T_1$. We also have
 \begin{center}
 	$
 	\begin{array} {rl}
 	z_X(S)	=	&	z_X({T_1})+2z_X(2^{E\backslash\{e,f_1,f_3\}}\backslash T_1)+z_X({T_2})+2z_X(2^{E\backslash\{e,f_1,f_3\}}\backslash T_2)+z_X(V)+z_X(W)		\\
 	=	&	(2^{n-1-|X|})-z_{X}(T_2)-z_X(T_1)+z_X(V)+z_{X}(W)\\
 	= &(2^{n-1-|X|})-z_{X}(T_2)+z_{X}(W).
 	\end{array}
 	$
 \end{center} 
 Since $z_X(W),z_X({T_2})\leq 2^{n-3-|X|}$ and $z_X(S)$ is a power of two, we have $z_X(W)=z_X(T_2)$ for all $X$, and so $W=T_2$. Hence $e$ and $f_1$ are deletable elements of $S$.

 Therefore every powerful set of order $n\geq 2$ and rank $n-1$ contains a deletable element. $\hfill\Box$
\\
 
\begin{pro} \label{pro:frame_existance}
	A powerful set $S$ of order $n\geq 2$ has rank $n-1$ if and only if it is contains a star.
\end{pro} 
\pf We prove this by induction on $n$.

The two powerful sets of order 2 and rank 1 each contain a star.

By Theorem \ref{thm:deletable_rank}, $S$ contains a deletable element $e$. If $e$ is not a coloop, then $S\backslash{e}=2^{E\backslash\{e\}}$, and so $S=S/e+\star_e$ by definition. If $e$ is a coloop, then $S\backslash{ e}$ has order $n-1$ and rank $n-2$, and by induction $S\backslash{e}$ contains a star element $d\in E\backslash\{e\}$, with $S\backslash{e}=T+\star_d$ for some powerful set $T$. But $S=S\backslash{ e}+\circ_e^*=(T+\star_d)+\circ_e^*$. It follows from the definitions of coloops and stars that for any powerful set $T$ we have $(T+\star_d)+\circ_e^*=(T+\circ_e^*)+\star_d$, so $d$ is a star of $S$. 

Finally, a powerful set $S$ of order $n$ that contains a star $e$ has rank $n-1$. $\hfill\Box$

This answers the conjecture from \cite{farr2017} as to whether every powerful set of order $n\geq 2$ and rank $n-1$ contains a star in the affirmative.

Proposition \ref{pro:frame_existance} gives a necessary and sufficient condition for a powerful set to contain a star, but does not indicate which element is a star. It follows from the definition of a star that an element $e\in E$ is a star of a powerful set $S\subseteq 2^E$ if and only if $r_S(X)=|X|$ for all $X\subseteq E\backslash\{e\}$. So it is possible to determine whether an element is a star with a number of rank evaluations exponential in the size of the ground set.  This is in contrast to Theorem \ref{thm:extension_rank}, where it can be determined in a constant number of rank evaluations whether or not an element is a loop, coloop or frame.

This same result gives a bijection between powerful sets of order $n$ and rank $n-1$ to powerful sets of order $n-1$, giving the following result.
\begin{cor}
 	The number of powerful sets of order $n\geq 2$ and rank $n-1$ is equal to the number of powerful sets of order $n-1$.
 \end{cor}

	\section{Characterising Non-linearity} 
	\label{sec:nonlinear}	
	The rank functions of linear powerful sets are precisely the binary matroid rank functions. All matroid rank functions satisfy the matroid rank axioms, one of which is \emph{subcardinality}. Namely, a function $r:2^E\to\mathbb{N}\cup\{0\}$ is \emph{subcardinal} if, for all $X\subseteq E$, $r(X)\leq |X|$. 
	
	The following theorem shows that, for powerful sets and multisets, only linear powerful sets have subcardinal rank functions.
	
	\begin{thm}\label{thm:linear}
		Let $S$ be a powerful multiset with ground set $E$. If $r_S(X)\leq |X|$ for all $X\subseteq E$, then $S$ is a linear powerful set.
	\end{thm}
	
	\pf
		We give a proof by contradiction. Assume there exists a nonlinear powerful multiset $S$ such that $r_S(X)\leq |X|$ for all $X\subseteq E$.  We take $S$ to be a powerful multiset of minimum order, i.e., every nonlinear powerful multiset of order less than $S$ does not have a non-subcardinal rank function.
		By Theorem \ref{thm:multiset_scale}, we may take $S$ to be such that the multiplicity of the empty set in $S$ is one.
	
	For a nonlinear powerful multiset $S$ over the ground set $E=\{e\}$ we must have $r_S(\{e\})\notin \{0,1\}$, and so $r_S(\{e\})> 1$. So no nonlinear powerful multiset of order 1 has a subcardinal rank function. We must have that the order of $S$ is at least two.

	Suppose there is some element $e$ which has rank 0 in $S$, $r_S(\{e\})=0$. 	By definition, $r_S(\{e\})=\log_2(\sum_{Y\subseteq E}f_S(Y)/\sum_{Y\subseteq E\backslash\{e\}}f_S(Y))=0$, so $\sum_{Y\subseteq E}f_S(Y)=\sum_{Y\subseteq E\backslash\{e\}}f_S(Y)$. Therefore, $f_S(Y\cup\{e\})=0$ for each $Y\subseteq E\backslash{\{e\}}$. For any $X\subseteq E\backslash\{e\}$, 

\[	\begin{array}{rl}
	    r_S(X\cup\{e\})= &  \log_2\left( \dfrac{\sum_{Y\subseteq E}f_S(Y)}{\sum_{Y\subseteq E\backslash{(X\cup\{e\})}}f_S(Y)}\right) \\	 
	     = &  \log_2\left( \dfrac{\sum_{Y\subseteq E}f_S(Y)}{\sum_{Y\subseteq E\backslash{X}}\left(f_S(Y)+f_S(Y\cup\{e\})\right)}\right) \\
	     = &  \log_2\left( \dfrac{\sum_{Y\subseteq E}f_S(Y)}{\sum_{Y\subseteq E\backslash{X}}f_S(Y)}\right)\\
	     =&r_S(X).\\
	\end{array}\]
	
	Hence, for any $X\subseteq E\backslash \{e\}$, we must have $r_S(X)=r_S(X\cup \{e\})$. Deleting $e$ gives $S\backslash e$, and for all $X\subseteq E\backslash \{e\}$ we have $r_{S\backslash e}(X)=r_S(X)\leq |X|$, as $r_S$ is subcardinal. Hence the rank function of  $S\backslash{e}$ is subcardinal, and as $S\backslash{e}$ has order less than $S$, $S\backslash e$ is a linear powerful set. 

	However, for any $X\subseteq E$, $r_S(X)=r_S(X\backslash\{e\})=r_{S\backslash e}(X\backslash\{e\})=r_{S\backslash e+\circ_e}(X)$, and so $S=S\backslash{e}+\circ_e$ is a linear powerful set (see \cite{farr2017}). As $S$ is nonlinear, we must then have that $r_S(\{e\})\neq 0$ for each $e\in E$. Hence $r_S(\{e\})=1$ for all $e\in E$.
	
	Contracting $e$ gives $S/e$, and for all $X\subseteq E\backslash e$ we have $r_{S/e}(X)=r_S(X\cup \{e\})-r_S(\{e\})\leq |X|$, as $r_S$ is subcardinal. Hence  $S/{e}$ has a subcardinal rank function and order less than $S$, so $S/e$ is a linear powerful set. As $S$ is nonlinear, $e$ is not a coloop of the linear powerful set $S/e$ (see \cite{farr2017}). 
	Since $S\backslash e$ is also a linear powerful set, all elements in $S\backslash e$ must appear with multiplicity one, and so $S$ is a powerful set. We now have that $S\subseteq 2^E$ must be a nonlinear powerful set with subcardinal rank function.
	
    In summary, we must have that $S\subseteq 2^E$ is a nonlinear powerful set with subcardinal rank function, and that deletion or contraction by any element gives a linear powerful set. In particular, we have that every element of $S$ is deletable. 
	
	Firstly, note that $S$ contains no sets of size one, for if $\{e\}\in S$, then by Theorem \ref{thm:coloops} $e$ is a coloop of $S$, and so $S=S/e+\circ^*_e$ would be linear. 
	
	As $S$ is nonlinear, there exist $X,Y\in S$ such that the symmetric difference $X\Delta Y\notin S$.

	Let $e\in E$. We have that $X\backslash\{e\}$ and $Y\backslash\{e\}$ in $S\backslash e$. As $S\backslash{ e}$ is linear, we have 
	\[(X\Delta Y)\backslash \{e\}=	(X\backslash \{e\}) \Delta (Y\backslash \{e\})  \in S\backslash e.	\]
	
Since $(X\Delta Y)\backslash \{e\} \in S\backslash e$, at least one of $X\Delta Y$ or $X\Delta Y\Delta \{e\}$ is a member of $S$. But $X\Delta Y \notin S$, so $X\Delta Y \Delta \{e\}\in S$.
	
Let $e'\in E\backslash\{e\}$. By our previous observations, we have $X\Delta Y \Delta \{e'\}\in S$, and so 
	
	\[\begin{array}{rl}
(X\backslash \{e\} )\Delta (Y\backslash \{e\} ) \Delta \{e'\}=	(X\Delta Y \Delta \{e'\})\backslash \{e\}&\in S\backslash e\\
	(X\Delta Y)\backslash \{e\} =(X\backslash \{e\}) \Delta (Y\backslash \{e\}) &\in S\backslash e\\
	((X\backslash \{e\} )\Delta (Y\backslash \{e\}) \Delta \{e'\})\Delta((X\backslash \{e\}) \Delta (Y\backslash \{e\}))=\{e'\} &\in S\backslash e\\
	\end{array}	\]

	As $S$ has no sets of size one, and $\{e'\}\in S\backslash e$, we must have $\{e,e'\}\in S$. This holds for any pair of elements in $E$, and so $S$ contains all subsets of $E$ of size two.  As $S$ contains no sets of  of size one the set of minimal nonempty members, or cocircuits, of $S$ is
	\[	\mathcal{C}(S)=\{X\subseteq E:|X|=2\}.	\]
	
	By Theorem \ref{thm:cocircuits}, $S$ must be isomorphic to the binary matroid $M(C_{|E|})$ (the cycle matroid of the $|E|$-vertex cycle). But this is a contradiction, as $S$ is not be linear. Hence no such powerful set exists. 
		
	$\hfill\Box$
	
As every linear powerful set has a subcardinal rank function, we have the following.
\begin{cor}
	A powerful set $S\subseteq 2^E$ is linear if and only if $r_S(X)\leq |X|$ for all $X\subseteq E$.
	\end{cor}
	
	\section{Conclusion and Further Work}

    In this paper we showed how degenerate elements (loops, coloops, frames and stars) relate to the rank function of a powerful set. We found that a fixed number of evaluations of the rank function can determine if an element is a loop, coloop or frame. We found that the rank of a powerful set determines if it contains a star.
    These results prove the conjectures from \cite{farr2017}.
	
	One avenue of investigation of the rank function of a powerful set would be to look at partial evaluations of a Tutte-Whitney type polynomial (see \cite{ellis2011}).  Results from \cite{farr93,farr04} generalise the characteristic polynomial, percolation probability, weight enumerator and Potts model partition function as partial evaluations of a binary function's Whitney rank generating function. In the case of powerful sets and powerful multisets these partial evaluations are polynomials. We can investigate properties of these polynomials, such as combinatorial interpretations of coefficients and locations of roots.

    Polymatroids are a natural generalisation of matroids (see, e.g. \cite[ \S 18.2]{welsh76}), and have a rank function which is not necessarily subcardinal. A powerful multiset can be constructed with the rank function from the binary vector space associated with a binary polymatroid. We found that the intersection of the classes of matroids and powerful multisets is the class of binary matroids. A natural problem would be to determine the powerful multisets with a polymatroid rank function, and determine whether they are solely the binary polymatroids.

    \subsection*{Acknowledgements} 
    The author would like to thank Graham Farr and Kerri Morgan for their supervision and help in developing this paper, and Keisuke Shiromoto for helpful discussions, in particular with Theorem \ref{thm:linear}.

\end{document}